\documentclass{article}

\usepackage[all]{xy}
\usepackage{amsfonts}
\usepackage{amssymb}
\usepackage{bbm}

\newtheorem{theorem}{Theorem}[section]

\newtheorem{definition}[theorem]{Definition}

\newtheorem{problem}[theorem]{Problem}

\newcommand{\xn}{\mathbb{N}}

\newcommand{\xq}{\mathbb{Q}}
\newcommand{\xf}{\mathbb{F}}

\newcommand{\xz}{\mathbb{Z}}

\newcommand{\mat}[2]{\mbox{Mat}_{#1}(#2)}

\title{Fast point counting on genus two curves \\
in characteristic three}

\author{Robert Carls}

\begin{document}

\maketitle

{\small 
\begin{center}
\begin{tabular}{l}
Robert Carls \\
\texttt{robert.carls@uni-ulm.de} \\
Universit\"{a}t Ulm \\
Institut f\"{u}r Reine Mathematik \\
D-89069 Ulm, Germany 
\end{tabular}
\end{center}
}

\begin{abstract}
\noindent
In this article we give the details of an
effective point counting algorithm for
genus two curves over finite fields of characteristic three.
The algorithm has an application in the context of curve based
cryptography.
One distinguished property of the algorithm is that its
complexity depends quasi-quadratically
on the degree of the finite base field. Our algorithm is
a modified version of
an earlier method that was developed in joint work with Lubicz.
We explain how one can alter the original algorithm, on the basis
of new theory, such that it can be used to efficiently count points
on genus two curves over large finite fields.
Examples of cryptographic size
have been computed using an experimental Magma
implementation of the algorithm which has been programmed by the
author.
Our computational results show that the quasi-quadratic
algorithm of Lubicz and the author, with
some improvements, is practical and relevant for cryptography. 
\end{abstract}

\section{Introduction}

\noindent
In this article we give the details of an
effective point counting algorithm for
genus two curves over finite fields of characteristic three, the
complexity of which depends quadratically on the degree of the finite
base field.
Our algorithm
is a modified version of
an earlier method that was developed in joint work with Lubicz \cite{cl09}.
The main purpose of this paper is to show that the original algorithm
of Lubicz and the author, with some improvements based on new theory,
is practical for genus $2$ curves over finite fields of
cryptographiy size. We conclude that our point counting
algorithm is relevant for
curve based cryptography. The importance
of genus $2$ cryptosystems comes from the fact that the size of the base field
can be chosen significantly smaller, namely half the size,
than in the case of elliptic curve systems at the same security level. 
This makes genus $2$ curves attractive for applications on crypto-devices
with limitations on the computing resources.
On the other hand, the quasi-quadratic dependency on the size of the
base field makes it possible to compute curves over huge finite fields
which are suitable for cryptography on the highest level of security.
\newline\indent
Our point counting algorithm can be used for the generation
of the key data that is necessary
for public key cryptography on the basis of low genus curves
over finite fields.
Usually, the key data of an algebraic curve cryptosystem consists
of the following objects
\begin{enumerate}
\item[{\rm (I)}]
a non-singular projective curve $C$ over a finite field
$\xf_q$ with $q$ elements,
\item[{\rm (II)}]
a computational model of the Jacobian group variety $J_C$ of the curve $C$,
\item[{\rm (II)}]
points $P,Q \in J_C(\xf_q)$ such that there exists
an $m \in \xn$ with $[m](P)=Q$.
\end{enumerate}
Given the above data one can for example encrypt or sign data using the
generic ElGamal method. For a detailed discussion of curve
based cryptography we refer to \cite{ko98}. 
The problem of finding the number $m$ from the given tuple
\begin{eqnarray}
\label{key}
(\xf_q,C,J_C,P,Q)
\end{eqnarray}
is called the \emph{discrete logarithm problem}.
A curve $C$ over a finite field $\xf_q$
as above is considered as secure if
the cardinality $q$ of the finite field $\xf_q$ is
such that the discrete logarithm problem
in the group of $\xf_q$-rational points $J(\xf_q)$
is computationally infeasible.
Provided that $q$ has been chosen suitably large,
giving the required security level, one has to make sure that
the number of $\xf_q$-rational points $\#J_C(\xf_q)$ of the
Jacobian has a large prime factor
of bit size almost equal to $g \cdot \log_2(q)$, where $g$ denotes the
genus of the curve $C$.
Under these assumptions, the generic methods for
solving the discrete logarithm problem are not
applicable.
The choice of a sufficiently large finite field $\xf_q$ is
a matter of finding a trade off between the desired security level and
the efficiency of en- and decryption functionality.
To check whether a given key data (\ref{key}) is secure
in the above sense
one has to compute and factorize the group
order $\#J_C(\xf_q)$.
In this article we discuss the following problem in a special case.
\begin{problem}
\label{pcproblem}
For a given finite field $\xf_q$ and a given curve $C$ over $\xf_q$,
compute the number $\#J_C(\xf_q)$.
\end{problem}
\noindent
Lubicz et al.~have proven in \cite{cl09} and \cite{ll06} that there exists
a quasi-quadratic algorithm which solves the Problem \ref{pcproblem} in the
case where the curve $C$ is ordinary and hyperelliptic.
Since their method depends polynomially on the characteristic
of the finite field, in practice it is limited to small
characteristics. In this article we describe some improvements of the
original algorithm, based on new theory, and an implementation of the improved
algorithm for ordinary genus $2$ curves over finite fields
of characteristic $3$.
This implementation has enabled us to compute examples of
cryptographic size in a reasonable amount of time. Our
computational results show that the method of Lubicz and the
author, with some modification, is practical and relevant for cryptography. 
\newline\indent
Let us now recall the precise result (compare \cite[Th.3.1]{cl09})
in the special case that is discussed in this article.
\begin{theorem}
\label{qq}
One can give an effective algorithm which computes
for an explicitly given non-singular ordinary genus $2$ curve $C$,
which is defined over a finite field $\xf_q$ of characteristic $3$,
the number $\#J_C(\xf_q)$ in time $O \big( \log(q)^{2+ \epsilon}
\big)$ for all $\epsilon >0$.
\end{theorem}
\noindent
We give the algorithm of Theorem \ref{qq} in Section \ref{algo}.
Due to a limited amount of space,
we don't give a complete proof of the correctness of our algorithm
in here.
\newline\indent
Let us remind the reader
that a distinguished property of our algorithm is given by the
fact that it is quasi-quadratic in the degree of the finite field.
Other algorithms by Kedlaya \cite{ke01}, Lauder and Wan
\cite{lw08} are just quasi-cubic,
which makes a difference in practice if the size of the finite field
is very large.
The algorithm of this article performs well for genus two curves
over finite fields of size far beyond the standards of state-of-the-art
hyperelliptic curve cryptography.

\subsection*{Leitfaden}

In Section \ref{algo} we give the details of the algorithm whose
existence is claimed in Theorem \ref{qq}.
We give examples that have been computed
using our algorithm in Section \ref{examples}.

\section{Algorithm}
\label{algo}

\noindent
In this section we give the algorithm which is subject to Theorem
\ref{qq}.
By $\xf_q$ we denote a finite field
with $q$ elements which is of characteristic $3$. All
computations in $\xf_q$ are supposed to be performed with polynomials
in $\xf_3[x]$ modulo a fixed irreducible monic polynomial
$\bar{f}$ over $\xf_3$ with $\mathrm{deg}(\bar{f})=\log_3(q)$ using
a fast polynomial arithmetic.
With $\xz_q$ we denote the ring $\xz_3[x]$ modulo the ideal
which is generated by a monic polynomial $f \in \xz[x]$ such that
$\bar{f} \equiv f \bmod 3$ and $\mathrm{deg}(f)=\log_3(p)$.
The ring $\xz_q$ is called the ring of Witt vectors with values
in $\xf_q$.
We say that we are given an element $x \in \xz_q$ with precision
$m$ if we have computed a bit string which represents
the truncated $3$-adic number $x$ modulo $3^m$.
Let $\sigma \in
\mathrm{End}_{\xz_3}(\xz_q)$ denote the unique lift
of the $3$-rd power Frobenius of $\xf_q$.
\newline\indent
The input of the point counting algorithm consists,
first, of a finite field $\xf_q$ in the above presentation
and, secondly, of an ordinary
curve $C$ which is given by an equation of the form
\begin{eqnarray}
\label{genus2}
y^2=x(x-1)(x-e_1)(x-e_2)(x-e_3)
\end{eqnarray}
where $e_1,e_2,e_3 \in \xf_q \setminus \{0,1\}$
are pairwise distinct.
We note that with only slight modification
our algorithm can also be applied to hyperelliptic
genus $2$ curves of a more general form.
The algorithm outputs the characteristic polynomial $\chi$
of the Frobenius endomorphism of the Jacobian variety $J_C$
of $C$ which is given by the $q$-th powering map.
From this one can obtain the number of $\xf_q$-rational
points $\#J_C(\xf_q)$ of $J_C$ by evaluating the polynomial
$\chi$ at the value $1$.
In the following sections we describe all steps of the
point counting algorithm in detail.
\begin{enumerate}
\item[{\rm (I)}]
Compute a $6$-theta null point $T_6$ of the curve $C$.
\item[{\rm (II)}]
Canonically lift the $6$-theta null point $T_6$
to a canonical theta null point $\tilde{T}_6$ with sufficiently
high precision.
\item[{\rm (III)}]
Compute the norm $\mathrm{Norm}( \delta )$ of the
determinant $\delta$ of a lift of the relative Verschiebung
in terms of the coefficients of $\tilde{T}_6$.
\item[{\rm (IV)}]
Reconstruct the characteristic polynomial $\chi$
from the approximated value for $\mathrm{Norm} ( \delta )$.
\end{enumerate}

\subsection{Computation of theta null points}

\noindent
In this section we explain how to compute the $6$-theta null point
of a curve given by an equation of the form (\ref{genus2}).
First, one computes a $2$-theta null point
\begin{eqnarray}
\label{theta2}
T_2=(b_{00},b_{01},b_{10},b_{11})
\end{eqnarray}
possibly over an extension field,
using the following
classical \emph{Thomae formulae}
\begin{eqnarray}
\label{thomae}
&& b_{00}=1 \\
\nonumber && b_{01}=\sqrt[4]{\frac{(e_1-e_4)(e_2-e_5)(e_3-e_4)}
{(e_1-e_5)(e_2-e_4)(e_3-e_5)}} \\
\nonumber && b_{10}=\sqrt[4]{\frac{(e_1-e_2)(e_1-e_4)}
{(e_1-e_3)(e_1-e_5)}} \\
\nonumber && b_{11}=\sqrt[4]{\frac{(e_1-e_2)(e_2-e_5)(e_3-e_4)}
{(e_1-e_3)(e_2-e_4)(e_3-e_5)}}
\end{eqnarray}
We note that for all possible roots in the formulas (\ref{thomae})
one gets a valid
$2$-theta null point.
In fact, the $2$-theta null point (\ref{theta2})
computed by means of the formulae (\ref{thomae})
belongs to an abelian variety which is $2$-isogenous to
the Jacobian of the curve $C$ defined by the equations (\ref{genus2}).
\newline\indent
Let us fix some notation. We denote $Z_n=(\xz / n \xz)^2$
for a natural number $n \geq 1$. Suppose further that we have chosen
embeddings $Z_n \hookrightarrow Z_m$, whenever $n|m$.
\newline\indent
Now once the $2$-theta null point is established, one needs to
extend the latter point to a smooth $6$-theta null point which we
denote by $T_6=(a_u)_{u \in Z_6}$.
The point $T_6$ lies in the zero locus of the following
equations
\begin{enumerate}
\item
\emph{symmetry relations}
\begin{eqnarray}
\label{symmetry}
Y_u=Y_{-u}, \quad \forall u \in Z_6
\end{eqnarray}
\item
\emph{Riemann relations}
\begin{eqnarray}
\label{riemann}
\lefteqn{\sum_{t \in Z_2} Y_{u_1+t} Y_{w_1+t} \cdot \sum_{t \in Z_2} Y_{z_1+t} Y_{y_1+t}}\\
\nonumber && = \sum_{t \in Z_2} Y_{u_2+t} Y_{w_2+t} \cdot \sum_{t \in Z_2} Y_{z_2+t}
Y_{y_2+t}
\end{eqnarray}
where $(u_i,w_i,z_i,y_i) \in Z_{6}^4$ for $i=1,2$ are equivalent quadruples.
\end{enumerate}
We consider
quadruples $(v_i,w_i,x_i,y_i) \in Z_{6}^4$ ($i=1,2$) as
equivalent if there exists a permutation matrix $P \in \mat{4}{\xz}$ such that
\begin{eqnarray*}
\lefteqn{(v_1+w_1,v_1-w_1,x_1+y_1,x_1-y_1)}\\
&&=(v_2+w_2,v_2-w_2,x_2+y_2,x_2-y_2)P
\end{eqnarray*}
The $2$-theta null point $T_2$ can be extended to a $6$-theta null
point $T_6$ using the above relations (\ref{symmetry}) and (\ref{riemann}).
We set
\[
a_{00}=b_{00}, \quad a_{03}=b_{01}, \quad a_{30}=b_{10}, \quad
a_{33}=b_{11}.
\]
Specializing the equations (\ref{symmetry}) and (\ref{riemann})
at $(a_{00},a_{03},a_{30},a_{33})$
we obtain a zero dimensional algebraic set in the variables
$\{ Y_u \}_{u \in Z_6 \setminus Z_2}$ (see \cite[Th.2.7]{cl09}).
The finiteness of this algebraic set enables us
to solve for a completed $6$-theta null point $T_6=(a_u)_{u \in Z_6}$.
\begin{theorem}
\label{existence}
There exists a smooth $6$-theta null point $T_6=(a_u)_{u \in Z_6}$ which
forms an extension of the $2$-torsion component
$(a_{00},a_{03},a_{30},a_{33})$.
\end{theorem}
\noindent
The proof of this fact involves sophisticated theory, so we
are not able to give it in here.
Here, smoothness means that some Jacobian criterion with respect
to the Riemann relations and additional correspondence
relations is satisfied at the point $T_6$. We will make that precise in Section \ref{canlift}.
We note that the smooth extension $T_6$ of Theorem
\ref{existence} belongs
to an abelian surface which is isogenous to the Jacobian of the curve $C$.
\newline\indent
The following method forms an important improvement
of the original algorithm as presented in \cite[$\S$3]{cl09}.
Instead of solving the relations (\ref{symmetry}) and (\ref{riemann})
for a smooth $6$-theta null
point $T_6$, we restrict to a smaller system
of equations which makes the problem feasible in practice.
Now consider the following subset of the Riemann equations
in the variables $Y_{10},Y_{13},Y_{20},Y_{23}$ which is given
by the four equations
\begin{eqnarray}
\label{smallsystem}
0 & = &
a_{00}^3Y_{20}+a_{00}^2a_{03}Y_{23}+a_{00}^2a_{30}Y_{10}+a_{00}^2a_{33}Y_{13}\\
\nonumber &&+a_{00}a_{03}^2Y_{20}+a_{00}a_{30}^2Y_{20}+a_{00}a_{33}^2Y_{20}+a_{03}^3Y_{23}
\\
\nonumber &&+a_{03}^2a_{30}Y_{10}+a_{03}^2a_{33}Y_{13}+a_{03}a_{30}^2Y_{23}+a_{03}a_{33}^2Y_{23}\\
\nonumber &&+a_{30}^3Y_{10}+a_{30}^2a_{33}Y_{13}+a_{30}a_{33}^2Y_{10}+a_{33}^3Y_{13}+2Y_{10}^4\\
\nonumber && +Y_{10}^2Y_{20}^2
+Y_{10}^2Y_{13}^2+Y_{10}^2Y_{23}^2+2Y_{20}^4+Y_{20}^2Y_{13}^2\\
\nonumber &&+Y_{20}^2Y_{23}^2+2Y_{13}^4+Y_{13}^2Y_{23}^2+2Y_{23}^4 \\
\nonumber 0 & = &
2a_{00}^2a_{33}Y_{13}+2a_{00}a_{03}a_{30}Y_{13}+2a_{00}a_{03}a_{33}Y_{10}\\
\nonumber &&+2a_{00}a_{30}a_{33}Y_{23}+2a_{00}a_{33}^2Y_{20}+2a_{03}^2a_{30}Y_{10}\\
\nonumber &&+2a_{03}a_{30}^2Y_{23}+2a_{03}a_{30}a_{33}Y_{20}+2Y_{10}^2Y_{23}^2\\
\nonumber &&+Y_{10}Y_{20}Y_{13}Y_{23}+2Y_{20}^2Y_{13}^2 \\
\nonumber 0 & = &
2a_{00}^2a_{30}Y_{10}+2a_{00}a_{03}a_{30}Y_{13}+2a_{00}a_{03}a_{33}Y_{10}\\
\nonumber &&+2a_{00}a_{30}^2Y_{20}+2a_{00}a_{30}a_{33}Y_{23}+2a_{03}^2a_{33}Y_{13}\\
\nonumber &&+2a_{03}a_{30}a_{33}Y_{20}+2a_{03}a_{33}^2Y_{23}+2Y_{10}^2Y_{20}^2\\
\nonumber &&+Y_{10}Y_{20}Y_{13}Y_{23}+2Y_{13}^2Y_{23}^2 \\
\nonumber 0 & = & 2a_{00}^2a_{03}Y_{23}+2a_{00}a_{03}^2Y_{20}+2a_{00}a_{03}a_{30}Y_{13}\\
\nonumber &&+2a_{00}a_{03}a_{33}Y_{10}+2a_{00}a_{30}a_{33}Y_{23}+2a_{03}a_{30}a_{33}Y_{20}\\
\nonumber &&+2a_{30}^2a_{33}Y_{13}+2a_{30}a_{33}^2Y_{10}+2Y_{10}^2Y_{13}^2\\
\nonumber &&+Y_{10}Y_{20}Y_{13}Y_{23}+2Y_{20}^2Y_{23}^2
\end{eqnarray}
where it is assumed
that the finite field elements $a_{00},a_{03},a_{30}$ and $a_{33}$ have already
been computed.
\begin{theorem}
The system of equations $\mathrm{(\ref{smallsystem})}$ defines a zero dimensional
algebraic set.
\end{theorem}
The system (\ref{smallsystem}) is readily solved by a standard Groebner basis algorithm on a normal
desktop computer over finite fields of cryptographic size.
For simplicity, we now assume that the given values $a_{00},a_{03},a_{30},a_{33}$
are defined over the field $\xf_q$.
\begin{theorem}
\label{degree}
A smooth $6$-theta null point $T_6=(a_u)_{u \in Z_6}$ is $L$-rational
over a field extension $L$ of $\xf_q$ such that $[L:\xf_q]$
divides $48$.
\end{theorem}
Let us remark that in most cases the degree of the field
extension is small. Our computations show that
for many examples it has degree lower or equal than $3$.
As a consequence of Theorem \ref{degree},
one can compute for increasing extension degree
the set $\mathcal{S}$ of four tuples $(a_{10},a_{13},a_{20},a_{23})$
that form a solution of the system (\ref{smallsystem}).
The homogeneity of the space of solutions of the
Riemann relations with respect to the action of
the automorphism group of the theta group implies that the
solution set $\mathcal{S}$
also contains the quadruples
\begin{eqnarray*}
(a_{14},a_{11},a_{22},a_{25})\\
(a_{32},a_{31},a_{02},a_{01})\\
(a_{12},a_{15},a_{24},a_{21})
\end{eqnarray*}
Thus, by forming all possible combinations of solutions in $S$, one obtains
as set of possible candidates for the $6$-theta null point
$T_6=(a_u)_{u \in Z_6}$. This completes our exposition of
the initial computations in Step (I) of our algorithm.
\newline\indent
Finally, let us give some further details of our implementation.
One can quickly test whether a candidate for $T_6$ is a valid $6$-theta
null point using the special theta relation
\begin{eqnarray*}
0 & = & 2a_{00}a_{10}a_{01}a_{31}+2a_{00}a_{20}a_{31}^2
        +2a_{00}a_{13}a_{02}a_{31}\\
&& +2a_{00}a_{23}a_{31}a_{32} +2 a_{03}a_{10}a_{01}a_{32}
   +2a_{03}a_{20}a_{31}a_{32} \\
&& +2a_{03}a_{13}a_{02}a_{32}+2a_{03}a_{23}a_{32}^2
   +2a_{30}a_{10}a_{01}^2 \\
&& +2a_{30}a_{20}a_{01}a_{31}+2a_{30}a_{13}a_{01}a_{02}
   +2a_{30}a_{23}a_{01}a_{32} \\
&& +2a_{33}a_{10}a_{01}a_{02}+2a_{33}a_{20}a_{02}a_{31}+2a_{33}a_{13}a_{02}^2 \\
&&
   +2a_{33}a_{23}a_{02}a_{32}+a_{10}^2a_{25}a_{21}+a_{10}a_{20}a_{11}a_{21}\\
&& +a_{10}a_{20}a_{25}a_{15}+a_{10}a_{13}a_{22}a_{21}+a_{10}a_{13}a_{25}a_{24} \\
&& +a_{10}a_{23}a_{14}a_{21}
   +a_{10}a_{23}a_{25}a_{12}+a_{20}^2a_{11}a_{15}\\
&&
   +a_{20}a_{13}a_{11}a_{24}+a_{20}a_{13}a_{22}a_{15}+a_{20}a_{23}a_{11}a_{12}\\
&&
   +a_{20}a_{23}a_{14}a_{15}+a_{13}^2a_{22}a_{24}+a_{13}a_{23}a_{14}a_{24}\\
&& +a_{13}a_{23}a_{22}a_{12}+a_{23}^2a_{14}a_{12}
\end{eqnarray*}
The smoothness of a candidate for $T_6$
is tested by computing the rank of the Jacobian matrix
with respect to the Riemann relations, taken together with
the correspondence relations (\ref{correq1}) and (\ref{correq2}) that
we introduce at a later point.
We will give a precise formulation of the smoothness criterion
in Section \ref{canlift}.
We remark that a different method for the computation of the $6$-theta null
point is suggested in \cite{fl09}.

\subsection{Canonical lifting}
\label{canlift}

\noindent
We use the notation of the preceding section.
The computation of the canonical lifted $6$-theta null point $T_6$
is realized by applying a Hensel lifting algorithm
to a system of equations that we define in the following.
\newline\indent
Consider the system of \emph{correspondence relations}
\begin{eqnarray}
\label{correq1}
\sum_{t \in Z_{6}, \, 3t=u} X_w Y_t
= \sum_{s \in Z_{6}, \, 3s=w} X_u Y_s
\end{eqnarray}
where $w,u \in Z_2$, and
\begin{eqnarray}
\label{correq2}
\lefteqn{\sum_{z \in Z_2} X_{x_1+z} X_{y_1+z}
\cdot \sum_{u \in Z_{6}}
X_{v_2+3u} Y_{w_2+u}}\\
\nonumber && = \sum_{z \in Z_2} X_{x_2+z} X_{y_2+z} \cdot \sum_{u \in Z_{6}}
X_{v_1+3u} Y_{w_1+u}
\end{eqnarray}
where $(x_i,y_i,v_i,w_i) \in S$ ($i=1,2$) and
$S$ is defined as the set of all $4$-tuples $(x,y,v,w) \in Z_{6}^4$ such
that the sets $\{ x+y, x-y \}$ and $\{ v+3w, v-3w \}$
are equal and contained in $Z_3$.
\newline\indent
By general theory (compare \cite{ca07} and \cite{ca05b})
there exists a canonical lift
$\tilde{T}_6=(\tilde{a}_u)_{u \in Z_6}$
of the $6$-theta null point $T_6=(a_u)_{u \in Z_6}$ to $\xz_q$.
\begin{theorem}
\label{frobtheta}
The points $\tilde{T}_6$ and
$\tilde{T}_6^{\sigma^2}=(\tilde{a}_u^{\sigma^2})_{u \in Z_6}$
satisfy the correspondence relations $\mathrm{(\ref{correq1})}$
and $\mathrm{(\ref{correq2})}$, if one evaluates
the variables $X_u$ and $Y_u$ with the values
$\tilde{a}_u$ and $\tilde{a}_u^{\sigma^2}$,
respectively.
\end{theorem}
For the subset of correspondence equations (\ref{correq2}) the
Theorem \ref{frobtheta} follows from
\cite[Th.2.1]{cl09}.
In the case of the equations (\ref{correq1}) a proof of the
Theorem \ref{frobtheta} can be found in the forthcoming preprint \cite{cm10}.
\newline\indent
Next we give a precise definition of the smoothness condition
that is subject to Theorem \ref{existence}.
It is convenient to 
use a short representation of $6$-theta null points
$(x_u)_{u \in Z_6}$ of the following shape
\begin{eqnarray*}
&& (x_{01},x_{02},x_{03},x_{10},x_{11},x_{12},x_{13},x_{14},
x_{15},\\
&& x_{20},x_{21},x_{22},x_{23},x_{24},x_{25},x_{30},
x_{31},x_{32},x_{33})
\end{eqnarray*}
which is justified by the symmetry equations (\ref{symmetry})
and the fact that
in almost all cases one can normalize
with respect to $x_{00}$.
We set
\begin{eqnarray*}
U=\{ 01,02,03,10,11,12,13,14,15,\\
20,21,22,23,24,25,30,
31,32,33 \}.
\end{eqnarray*}
By evaluating $Y_{00}$ at $1$ and by replacing, if necessary, $Y_u$ by
$Y_{-u}$ we can assume that the Riemann relations (\ref{riemann})
are given by a set $\mathcal{R}$ of polynomials
in the variables $Y_u$ where $u \in U$.
By the same procedure, and by evaluating $X_{00}$ at $1$, we obtain
from the correspondence equations (\ref{correq1}) and (\ref{correq2})
as set of polynomials $\mathcal{C}$ in the variables
$X_u$ and $Y_u$ where $u \in U$.
\begin{definition}
\label{smoothness}
We call a $\xf_q$-rational simultaneous zero
$(a_u)_{u \in U}$ of the polynomials in the set
$\mathcal{R}$ a smooth point,
if there exist polynomials
$f_1, \ldots,f_{19} \in \mathcal{R} \cup \mathcal{C}$ such that the matrix
of partial derivatives
\[
D_Y=\left(\frac{\partial f_i}
{\partial Y_u}\right), \quad i=1,\ldots,19 
\]
has non-zero determinant at the point $(a_u) \times (a^{9}_u)$,
where the index $u$ ranges over $U$.
\end{definition}
\noindent
It is straight forward to test computationally,
whether an $\xf_q$-rational
solution $(a_u)_{u \in Z_6}$ of the relations
(\ref{symmetry}) and (\ref{riemann}) is smooth in
the sense of Definition \ref{smoothness}. For example,
one can form the Jacobian matrix of all relations in
the set $\mathcal{R} \cup \mathcal{C}$ with respect to the
the variables $\{ Y_u \}$ and test whether the rank of
the resulting matrix is equal to $19$ at the
point $(a_u) \times (a^{9}_u)$.
\newline\indent
Now assume that we are given a smooth $6$-theta null point
$T_6=(a_u)_{u \in U}$.
To find polynomials $f_1, \ldots,f_{19}$ in $\mathcal{R}
\cup \mathcal{C}$ as in Definition \ref{smoothness},
one searches over all polynomials in $\mathcal{R}$ until
one has found $16$ relations such that their Jacobian
matrix has rank equal to $16$. Then one has to
find $3$ additional polynomials in $\mathcal{C}$
such that the vertical join of the Jacobian matrices
has rank $19$ in total.
As in Definition \ref{smoothness} we denote the
Jacobian matrix of the resulting polynomials
$f_1, \ldots, f_{19}$ with respect to the variables
$\{ Y_u \}_{u \in U}$ by $D_Y$.
The matrix of partial derivatives of these polynomials with respect
to the variables $\{ X_u \}_{u \in U}$
is denoted by $D_X$. We note that necessarily the
determinant of $D_X$ at $(a_u) \times (a^{9}_u)$
equals zero.
\newline\indent
We define a function $\Phi: \xz_q^{19} \times
\xz_q^{19} \rightarrow \xz_q^{19}$ by setting
\begin{eqnarray}
\label{liftfunction}
\Phi ( x,y )
=\big( f_1(x,y), \ldots, f_{19}(x,y) \big).
\end{eqnarray}
for all $(x,y)=(x_u)_{u \in U} \times (y_u)_{u \in U}
\in \xz_q^{19} \times \xz_q^{19}$. 
Suppose that we want to compute the canonical
lift $\tilde{T}_6=(\tilde{a}_u)_{u \in U}$
of the $6$-theta null point $T_6$
with given precision $m$.
Assume that we are given $\tilde{T}_6$
with precision $\lceil m/2 \rceil$.
By Theorem \ref{frobtheta} we have
\begin{eqnarray}
\label{nulleq}
\Phi \left( \tilde{T}_6, \tilde{T}^{\sigma^2}_6 \right) \equiv 0
\bmod 3^{\lceil m/2 \rceil}
\end{eqnarray}
Using Taylor expansion it follows from the
congruence (\ref{nulleq}) that
\begin{eqnarray}
\label{newnulleq}
0 \equiv \Phi \left( \tilde{T}_6 + 3^{\lceil m/2 \rceil} \cdot \Delta ,
\tilde{T}^{\sigma^2}_6 + 3^{\lceil m/2 \rceil} \cdot \Delta^{\sigma^2} \right)
\bmod 3^m
\end{eqnarray}
where $\Delta \in \xz_q^{19}$,
is equivalent to the congruence
\begin{eqnarray}
\label{artinschreier}
\lefteqn{ 0 \equiv \frac{1}{3^{\lceil m/2 \rceil}} \cdot
D_Y (\tilde{T}_6, \tilde{T}^{\sigma^2}_6)^{-1} \cdot \Phi \left( \tilde{T}_6, \tilde{T}^{\sigma^2}_6 \right) } \\
\nonumber && + D_Y( \tilde{T}_6, \tilde{T}^{\sigma^2}_6 )^{-1} \cdot D_X( \tilde{T}_6, \tilde{T}^{\sigma^2}_6 ) \cdot \Delta \\
\nonumber && + \Delta^{\sigma^2} \bmod 3^{\lceil m/2 \rceil}. 
\end{eqnarray}
Here we use the fact that the point $T_6=(a_u)_{u \in U}$, which is the
reduction of $\tilde{T}_6=(\tilde{a}_u)_{u \in U}$ modulo $3$,
is a smooth point, and consequently, the matrix
$D_Y( \tilde{T}_6, \tilde{T}^{\sigma^2}_6 )$
is invertible modulo $3^{\lceil m/2 \rceil}$.
Hence, by solving the generalized Artin-Schreier equation
(\ref{artinschreier}) one can compute a $\Delta \in \xz_q^{19}$ with
precision $\lceil m/2 \rceil$ which
solves the congruence (\ref{newnulleq}).
\newline\indent
In the following we describe an algorithm for the solution of
the above special type of generalized Artin-Schreier equation.
Again this is done by a Hensel lifting process.
Suppose that we are given a solution $\Delta \in \xz_q^{19}$ of the congruence
\[
\Delta^{\sigma^2}+A \cdot \Delta + v  \equiv 0 \bmod 3^{\lceil n/2 \rceil} 
\]
where $\Delta,v \in \xz_q^{19}$ and $A \in \mathrm{Mat} ( 19, \xz_q)$
is a square matrix which is singular modulo $3$.
The above congruence implies that
solving the congruence
\begin{eqnarray}
\label{as}
(\Delta+3^{\lceil n/2 \rceil} \cdot \epsilon  )^{\sigma^2}
+A \cdot (\Delta +3^{\lceil n/2 \rceil} \cdot \epsilon    ) + v  \equiv 0 \bmod 3^n 
\end{eqnarray}
where $\epsilon \in \xz_q^{19}$,
is equivalent to solving the congruence
\begin{eqnarray}
\epsilon^{\sigma^2}+A \cdot \epsilon + w \equiv 0 \bmod 3^{\lceil n/2 \rceil}
\end{eqnarray}
where
\[
w=\frac{1}{3^{\lceil n/2 \rceil}}
\cdot \big( \Delta^{\sigma^2}+A \cdot \Delta + v \big).
\]
The above calculations can be summarized in a lifting algorithm
for $6$-theta null points which
is based on the fact that it is computationally straight forward
to solve Artin-Schreier equations modulo $3$.
Using the above Hensel lifting principle one can compute
the canonical theta null point $\tilde{T}_6$ to given
precision in time depending quasi-linearly on the precision and the value
$\log_3(q)$.
We will specify the precision that we use in our point
counting algorithm in Section \ref{charpoly}.
\newline\indent
We omit a detailed description of the
method that we use to solve
an Artin-Schreier equation of the form
\begin{eqnarray}
\label{asfp}
\epsilon^{p^2}+\bar{A} \cdot \epsilon + \bar{w} \equiv 0 \bmod 3. 
\end{eqnarray}
where $\bar{A}$ is a singular matrix modulo $3$.
The solution of the congruence (\ref{asfp}) comes down
to solving a linear system modulo $3$.
Since it is straight forward to adapt the method described in
\cite[Algo.5.2]{ll06} to our situation, we don't give the details in here. 
This completes our description of the approximation of the
canonically lifted $6$-theta null point $\tilde{T}_6$
that is the main objective of Step (II) of our algorithm.

\subsection{Recovery of the characteristic polynomial}
\label{charpoly}

\noindent
We use the notation of the preceding sections.
For the rest of this
section let $\xf_q$ denote the field of definition of the
$6$-theta null point $T_6=(a_u)_{u \in Z_6}$.
Assume that we are given the canonically
lifted $6$-theta null point $\tilde{T}_6=(\tilde{a}_u)_{u \in Z_6}$ with
precision $m$.
Suppose that we have normalized $\tilde{T}_6$ such that
$\tilde{a}_{00}=1$.
Let $\pi_1,\pi_2$ be the $3$-adically invertible eigenvalues of the
absolute $q$-Frobenius endomorphism on the Jacobian variety
$J_C$ of the curve $C$ which is given by the equation (\ref{genus2}).
We set
\begin{eqnarray}
\label{conjugates}
\bar{\pi}_1=\frac{q}{\pi_1} \quad \mbox{and} \quad
\bar{\pi}_2=\frac{q}{\pi_2}.
\end{eqnarray}
Then the characteristic polynomial of Frobenius
is given by the following polynomial with $\xq$-coefficients
\[
\chi(T)=(T-\pi_1)(T-\pi_2)(T-\bar{\pi}_1)
(T-\bar{\pi}_2).
\]
We note that the product $\pi_1 \pi_2$ of eigenvalues
can be regarded as an element in $\xz_3$ in an obvious way.
We set
\[
\delta=1+2 ( \tilde{a}_{02}+ \tilde{a}_{20}+\tilde{a}_{22}+\tilde{a}_{24})^{\sigma^2}.
\]
The number $\delta$ is called the determinant
of relative Verschiebung.
\begin{theorem}
\label{nform}
One has
\begin{eqnarray*}
\mathrm{Norm}_{\xz_q/\xz_3} (\delta)
=\pm \pi_1 \pi_2
\end{eqnarray*}
\end{theorem}
An equivalent formula has been established
in \cite[Th.2.8]{cl09}. A purely algebraic proof
of Theorem \ref{nform} is given in the forthcoming preprint \cite{cm10}.
Theorem \ref{nform} implies that we can compute the product
of eigenvalues $\pi_1 \pi_2$ up to sign with given
precision.
This concludes our remarks regarding Step (III) of our algorithm.
\newline\indent
In the following we describe how one can compute a list of
candidates for the characteristic polynomial
$\chi(T)$, which is part of Step (IV) of our algorithm.
We note that the number of $\xf_q$-rational
points of the Jacobian variety $J_C$ of $C$ is given
by $\chi(1)$.
One can eliminate the false candidates for $\chi(T)$ by evaluating
at $1$ and performing point multiplications with random
points in the group $J_C(\xf_q)$. We remark that there is a well-known
algorithm for the addition of divisor classes in the
group $J_C(\xf_q)$. This is folklore, so we don't give the
details here.
In the following we ignore the field
extension that is necessary to compute a rational smooth $6$-theta
null point $T_6$. An extension of the base field of the curve $C$
can be compensated by taking appropriate roots of the eigenvalues
$\pi_1$ and $\pi_2$.
\newline\indent
Now let us briefly describe how one can compute the
characteristic polynomial $\chi(T)$ from the approximated
product of eigenvalues $\pi_1 \pi_2$.
Assume that we are given a $3$-adic number $\pi$
such that $\pi \equiv \pm \pi_1 \pi_2 \bmod 3^m$,
where the precision is chosen such that $m=2 \log_3(q)+2$.
The polynomial
\[
P_{\mathrm{sym}}(T)= ( T-\pi_1\pi_1+\bar{\pi}_1 \bar{\pi}_2 )
( T-\pi_1 \bar{\pi}_2 + \bar{\pi}_1 \pi_2)
\]
is called the symmetric polynomial associated to $\chi(T)$.
In order to compute the characteristic polynomial
$\chi(T)$, one first computes candidates for the symmetric polynomial
$P_{\mathrm{sym}}(T)$ in terms of $\pi$.
By the above discussion one has
\begin{eqnarray}
\label{sympoly}
P_{\mathrm{sym}}(T)=T^2-sT+qt
\end{eqnarray}
for some integers $s$ and $t$, whose absolute value
is smaller or equal to $9q$.
There exists an $s_0 \in \xz$, whose residue $\bar{s}_0$ modulo $9$ lies in
the interval $[ 0, \ldots, 8 ]$, such that
$s \equiv \pm \pi + \bar{s}_0q \bmod 9q$.
The algorithm for computing $s$ simply tries all of the above
possibilities for the residue $\bar{s}_0$ of $s_0$. For each
possible $\bar{s}_0$ one gets a
corresponding $s$, in terms of which we claim that one can compute the
parameter $t$. Since $|s|\leq 9q$, one can for every possible
integer $s$ compute
an exact value for $s_0$ by 
the formula $s_0=\frac{s-(\pi+\frac{q^2}{\pi})}{q}$.
Finally, one chooses $t \equiv \pi \cdot s_0 \bmod 9q$.
The above described procedure determines a list of integer pairs $(s,t)$
which give possible candidates for the
polynomial $P_{\mathrm{sym}}(T)$.
\newline\indent
Now assume that we are given roots $\alpha$
and $\beta$ of a candidate for the polynomial $P_{sym}(T)$ in a
suitable number field.
Let $\tau_1, \ldots, \tau_4$ denote the roots
of the polynomials
$P_1(T)=T^2-\alpha T+q^2$ and $P_2(T)=T^2-\beta T+q^2$,
in a suitable extension field of the rational numbers.
Then candidates for the values $\pm \pi_1^2$ and $\pm \pi_2^2$
can be computed up to sign as products $\tau_j \tau_k$, where
$j,k \in \{ 1, \ldots, 4 \}$.
By taking square roots one obtains candidates for the
eigenvalues $\pi_1$ and $\pi_2$.
The latter values determine the characteristic polynomial
$\chi(T)$ by the formulae (\ref{conjugates}).
This finishes the exposition of our point counting algorithm.

\section{Practical results}
\label{examples}

\noindent
In this section we give an
example of cryptographic size that was computed using our
algorithm.
A complete
documentation of the example is available on the author's
website \cite{caxx}.
Let $f(T)= T^{120} + T^4 + 2 \in \xf_3[T]$.
We denote by $\bar{T}$ the congruence class of the polynomial
$T$ modulo the modulus $f$.
Consider the hyperelliptic genus $2$ curve $C$
over the finite field $\xf_{3^{120}}=\xf_q[T]/(f)$ with defining equation
\newpage
\begin{eqnarray*}
\lefteqn{y^2 = x^5 + (2\bar{T}^{119} + \bar{T}^{116} +
  \bar{T}^{115}+\bar{T}^{114} + 2\bar{T}^{112}} \\
&& + 2\bar{T}^{109} + 2\bar{T}^{107} + \bar{T}^{104} + \bar{T}^{103} +
\bar{T}^{102} + \bar{T}^{101} \\
&& + 2\bar{T}^{96} + \bar{T}^{95} +2\bar{T}^{91}+ 2\bar{T}^{90} + 2\bar{T}^{88} +
\bar{T}^{87} \\
&& +2\bar{T}^{85} + 2\bar{T}^{84} + \bar{T}^{83} + \bar{T}^{81} +
\bar{T}^{80} + 2\bar{T}^{79} \\
&& + 2\bar{T}^{78} +
  \bar{T}^{77} + 2\bar{T}^{73} + 2\bar{T}^{71} + 2\bar{T}^{68} +
  2\bar{T}^{67}\\
&& + \bar{T}^{65} + 2\bar{T}^{63} + 
    \bar{T}^{62} + \bar{T}^{59} + 2\bar{T}^{57} + 2\bar{T}^{56} \\
&& + 2\bar{T}^{54} + \bar{T}^{53} + \bar{T}^{52} + 
    2\bar{T}^{51} + \bar{T}^{48} + 2\bar{T}^{47} \\
&& + \bar{T}^{46} + 2\bar{T}^{45} + 2\bar{T}^{43} + \bar{T}^{41} + 
    2\bar{T}^{40} + 2\bar{T}^{38} \\
&& + \bar{T}^{36} + 2\bar{T}^{35} + 2\bar{T}^{34} + 2\bar{T}^{32} + 2\bar{T}^{31} + 
    2\bar{T}^{30} \\
&& + \bar{T}^{29} + \bar{T}^{25} + 2\bar{T}^{24} + 2\bar{T}^{23} +
    \bar{T}^{22} + \bar{T}^{21} \\
&& + \bar{T}^{20} + \bar{T}^{19} + 2\bar{T}^{17} + \bar{T}^{16} +
    2\bar{T}^{15} + 2\bar{T}^{13} \\
&& + \bar{T}^{10} + 2\bar{T}^9 + 
    2\bar{T}^8 + \bar{T}^7 + \bar{T}^6 + 2\bar{T}^2 + 2\bar{T})x^4 \\
&& + (2\bar{T}^{119} + 2\bar{T}^{117} + 
    2\bar{T}^{116} + \bar{T}^{115}  + \bar{T}^{114} \\
&& + \bar{T}^{112} + 2\bar{T}^{111} + \bar{T}^{110} + 2\bar{T}^{106} + 
    \bar{T}^{105} + 2\bar{T}^{104} \\
&& + 2\bar{T}^{103} + \bar{T}^{102} + \bar{T}^{97} + 2\bar{T}^{96} + 2\bar{T}^{94} + 
    2\bar{T}^{93} \\
&& + \bar{T}^{91} + 2\bar{T}^{90} + \bar{T}^{89} + \bar{T}^{88} +
    2\bar{T}^{85} + \bar{T}^{83} \\
&&    + \bar{T}^{82} + 2\bar{T}^{81} + 2\bar{T}^{80} + 2\bar{T}^{78} +
    \bar{T}^{75} + \bar{T}^{74} \\
&& + 2\bar{T}^{71}
 + 2\bar{T}^{70} +\bar{T}^{67} + 2\bar{T}^{66} + 2\bar{T}^{65} +
 2\bar{T}^{64} \\
&& + \bar{T}^{63} + \bar{T}^{60}
   + \bar{T}^{59} + \bar{T}^{58} + 2\bar{T}^{56} + \bar{T}^{55} \\
&& + \bar{T}^{54} + 2\bar{T}^{51} + 2\bar{T}^{50}
 + \bar{T}^{49} + \bar{T}^{47} + \bar{T}^{46} \\
&& + 2\bar{T}^{44} + 2\bar{T}^{42} + 2\bar{T}^{39} + \bar{T}^{36}
   + 2\bar{T}^{33} + 2\bar{T}^{31} \\
&& + \bar{T}^{29} + \bar{T}^{28} + 2\bar{T}^{26} + \bar{T}^{25} + \bar{T}^{24}
    + \bar{T}^{23} + 2\bar{T}^{22} \\
&& + \bar{T}^{21} + \bar{T}^{19} + \bar{T}^{17} + 2\bar{T}^{16} + 2\bar{T}^{15}
 + 2\bar{T}^{12} \\
&& + 2\bar{T}^{11} + \bar{T}^9 + 2\bar{T}^7 + \bar{T}^5 + \bar{T}^4 +
 \bar{T}^3 \\
&& + \bar{T}^2 + 2\bar{T})x^3 + (2\bar{T}^{119} + \bar{T}^{118} +
 \bar{T}^{117} + 2\bar{T}^{115} \\
&& + 2\bar{T}^{114} + \bar{T}^{111} + 2\bar{T}^{108} + \bar{T}^{107} +
 \bar{T}^{105} + 2\bar{T}^{104} \\
&&+ 2\bar{T}^{103} +
 2\bar{T}^{101} + 2\bar{T}^{99} + 2\bar{T}^{98} + 2\bar{T}^{97} +
 \bar{T}^{96} \\
&& + 2\bar{T}^{94} + 2\bar{T}^{86}
 + \bar{T}^{84} + 2\bar{T}^{83} + \bar{T}^{82} + 2\bar{T}^{80} \\
&& +\bar{T}^{78} + \bar{T}^{77} + 2\bar{T}^{76}
 + 2\bar{T}^{75} + 2\bar{T}^{73} + \bar{T}^{72} \\
&& + 2\bar{T}^{71} + \bar{T}^{69} + 2\bar{T}^{68} + \bar{T}^{67}
 + \bar{T}^{65} + \bar{T}^{64} + \bar{T}^{62} \\
&& + 2\bar{T}^{61} + 2\bar{T}^{60} + \bar{T}^{59} + 2\bar{T}^{58}
 + 2\bar{T}^{55} + \bar{T}^{51} \\
&& + \bar{T}^{50} + 2\bar{T}^{49} + \bar{T}^{48} + 2\bar{T}^{47} + \bar{T}^{41}
 + \bar{T}^{40} \\
&& + 2\bar{T}^{39} + 2\bar{T}^{38} + 2\bar{T}^{37} + \bar{T}^{36} +
 2\bar{T}^{30} + \bar{T}^{28} \\
&& + 2\bar{T}^{27} + 2\bar{T}^{26} + 2\bar{T}^{24} + 2\bar{T}^{23} +
 2\bar{T}^{22} + 2\bar{T}^{21} \\
&&+ 2\bar{T}^{20}
 + \bar{T}^{17} + 2\bar{T}^{15} + \bar{T}^{14} + 2\bar{T}^{12} +
 2\bar{T}^{11} \\
&&+ 2\bar{T}^{10} + \bar{T}^8
 + 2\bar{T}^7 + \bar{T}^5 + 2\bar{T}^4 + 2\bar{T}^2 + 2)x^2 \\
&& + (2\bar{T}^{118} + 2\bar{T}^{115}
 + 2\bar{T}^{114} + 2\bar{T}^{110} + \bar{T}^{109} \\
&&+ \bar{T}^{108} + \bar{T}^{106} + \bar{T}^{105} 
    + \bar{T}^{104} + \bar{T}^{103} + \bar{T}^{102} \\
&& + \bar{T}^{99} + \bar{T}^{98} + \bar{T}^{96} + 2\bar{T}^{95} + 
    2\bar{T}^{94} + \bar{T}^{93} \\
&& + 2\bar{T}^{90} + 2\bar{T}^{89} + 2\bar{T}^{87} + \bar{T}^{86} + 2\bar{T}^{85} + 
    2\bar{T}^{83} \\
&& + \bar{T}^{82} + \bar{T}^{80} + \bar{T}^{79} + \bar{T}^{78} +
    \bar{T}^{77} + \bar{T}^{76} \\
&& + 2\bar{T}^{74} + 
    2\bar{T}^{73} + 2\bar{T}^{72} + \bar{T}^{70} + 2\bar{T}^{69} +
    2\bar{T}^{68}\\
&& + 2\bar{T}^{67} + \bar{T}^{66} + 
    2\bar{T}^{65} + \bar{T}^{62} + \bar{T}^{61} + \bar{T}^{57} \\
&& + 2\bar{T}^{56} + 2\bar{T}^{53} + 2\bar{T}^{52} + 
    \bar{T}^{51} + \bar{T}^{48} + \bar{T}^{47} \\
&& + \bar{T}^{46} + \bar{T}^{45} + \bar{T}^{44} + \bar{T}^{43}
  + \bar{T}^{42} + \bar{T}^{41} \\
&& + 2\bar{T}^{39} + 2\bar{T}^{38} + \bar{T}^{37} + \bar{T}^{35} + \bar{T}^{34}
  + \bar{T}^{33} \\
&& + \bar{T}^{32} + 2\bar{T}^{31} + 2\bar{T}^{30} + \bar{T}^{29} +
  \bar{T}^{28} + \bar{T}^{27} \\
&&  + 2\bar{T}^{26} + \bar{T}^{25} + \bar{T}^{24} + \bar{T}^{23} +
  \bar{T}^{22} + 2\bar{T}^{21} \\
&& + \bar{T}^{19} + 2\bar{T}^{17} + 2\bar{T}^{14} + \bar{T}^{13} +
  2\bar{T}^{12} + 2\bar{T}^{11} \\
&& + \bar{T}^7 + 2\bar{T}^6 +
  \bar{T}^5 + 2\bar{T}^3 + \bar{T}^2 + 2\bar{T})x 
\end{eqnarray*}
The number of $\xf_q$-rational points on the Jacobian variety $J_C$
of $C$ equals
\begin{eqnarray*}
&& 32292460179985540075152248365\\
&& 95391097003917060756603284118\\
&& 54046812502670061472170389646\\
&& 4902240351775536748901686160
\end{eqnarray*}
The group order $\#J_C(\xf_q)$ has a large prime factor of size
$369$ bits.
Also, by computing the minimal polynomials of the Igusa
invariants of the curve $C$, one can verify that $\xf_q$
is a minimal field of definition for the curve $C$.
Thus, the curve satisfies the requirements for
a cryptographically secure genus $2$ curve.
The computation of the group order, using our algorithm,
took $1394$ seconds (CPU time) on an Intel Core $2$ E7700
with $8$Gb memory.
Comparing the running time of our experimental implementation to
the built in Magma implementation of Kedlaya's algorithm for genus $2$
curves , one can see that our results
are reasonable.

\section{Summary and perspectives}

\noindent
In this article, we have given the details of
an effective quasi-quadratic algorithm for point counting
on ordinary genus $2$ hyperelliptic curves over finite fields of characteristic
$3$, which performs very well in practice.
Further improvement may be achieved
regarding the following open problems of theoretical nature.
\begin{enumerate}
\item
Can one modify the general algorithm given in \cite{cl09} such that
its complexity depends only polynomially on the
logarithm of the characteristic of the finite field ? 
Note that the original algorithm depends polynomially on $p^g$,
where $p$ is the characteristic and $g$ is the genus of the curve.
\item
Can one significantly reduce the number of variables in the canonical lifting
algorithm which is described in Section \ref{algo} ?
Some results that might turn out to be useful in this context
are documented in \cite{lr10}.
\item
By introducing coarse invariants, which are supposed to be
expressions in certain theta constants, can one avoid
the field extensions that in some cases are necessary to
obtain a rational theta null point ?
\end{enumerate}

\bibliographystyle{plain}

\end{document}